\documentclass[11pt]{amsart}

\usepackage{amsmath}
\usepackage{amssymb}
\usepackage{graphicx}
\usepackage[initials]{amsrefs}
\usepackage{fancyhdr}
\usepackage{subcaption}
\usepackage{algorithm}
\usepackage{algpseudocode}
\usepackage{tikz}
\usepackage{multirow}
\usetikzlibrary{shapes.geometric, arrows.meta, positioning}

\newtheorem{theorem}{Theorem}[section]

\theoremstyle{definition}
\newtheorem{definition}[theorem]{Definition}

\newtheorem{remark}[theorem]{Remark}

\numberwithin{equation}{section}

\title[RO Framework for priority-aware coverage]{Robust priority-aware coverage optimization for aerial sensor networks}

\author[V. Datta]{Vanshika Datta}
\address[V. Datta]{Department of Mathematics, Indian Institute of Technology Kharagpur, Kharagpur, West Bengal, India}
\email{{\tt vanshikadutta28@gmail.com}}

\author[C. Nahak]{C. Nahak}
\address[C. Nahak]{Department of Mathematics, Indian Institute of Technology Kharagpur, Kharagpur, West Bengal, India}
\email{\tt cnahak@maths.iitkgp.ac.in}

\author[J. C. Yao]{J. C. Yao$^*$}
\address[J. C. Yao]{Research Center for Interneural Computing, China Medical University Hospital, China Medical University, Taichung 40402, Taiwan}
\email{\tt yaojc@mail.cmu.edu.tw}

\keywords{Robust optimization; Radius of robust feasibility; Voronoi diagram; Directional sensor network; Data uncertainty\\
$^*$Corresponding author}

\begin{document}

\begin{abstract}
This article presents a priority-aware robust coverage optimization framework for an aerial sensor network under sensor location uncertainty. Each region is assigned a priority weight, and the objective is to maximize the weighted coverage while maintaining robustness against positional perturbations. A mathematical optimization model is developed by incorporating surveillance constraints and an RRF-based robustness formulation into the proposed framework. An efficient priority-aware robust orientation optimization (PAROO) algorithm is then proposed to determine the sensor orientations that maximize the weighted coverage objective. Experimental results on an airport-inspired surveillance scenario demonstrate that the proposed framework effectively directs sensing resources toward high-priority regions and achieves higher weighted coverage than representative baseline approaches, highlighting its practical applicability in security-sensitive environments.
\end{abstract}

\maketitle


\section{Introduction}
Wireless sensor networks (WSNs) have become an indispensable technology for monitoring and surveillance applications owing to their capability to continuously collect and transmit information from geographically distributed locations \cite{err_1, app_1}. Among various quality of service (QoS) metrics, coverage is one of the most fundamental, aiming to maximize the monitored region while efficiently utilizing sensing resources \cite{cvg_1}. Compared with omnidirectional sensors, directional sensor networks (DSNs) pose greater challenges due to their limited field of view, whereas aerial directional sensor networks (ASNs) further increase the complexity by introducing additional deployment flexibility and orientation constraints \cite{vs_1, vs_2}. Consequently, coverage optimization in ASNs has emerged as an important research topic with applications in surveillance, environmental monitoring, and critical infrastructure protection \cite{asn_1, asn_2}.

In practice, sensors are often deployed in harsh environments where obstacle occlusion, irregular boundaries, localization errors, environmental disturbances, and sensor failures may significantly degrade coverage quality \cite{err_1, err_2}. To improve network reliability under such uncertainties, numerous probabilistic, stochastic, fuzzy, and robust optimization approaches have been proposed \cite{u_1, u_2}. In parallel, several coverage optimization methods have demonstrated significant improvements in directional sensor coverage, see \cite{dsn_1, dsn_2}; however, these methods generally assume deterministic sensing environments. Conversely, existing robust optimization approaches typically treat every location within the monitored region as equally important \cite{vd, ro_1}.

Conventional coverage optimization assumes uniform importance across the region of interest. However, practical surveillance environments, such as airports, military installations, and industrial facilities, contain regions with substantially different monitoring requirements. For example, airport runways, terminals, fuel storage facilities, and air traffic control towers require significantly higher levels of surveillance than parking areas or open spaces. Recently, priority-aware coverage optimization has gained increasing attention by assigning higher importance to mission-critical regions, thereby improving surveillance effectiveness \cite{pa_1,pa_2}. Nevertheless, these studies generally neglect deployment uncertainty, while uncertainty-aware methods rarely incorporate regional priorities into their formulations. Consequently, the integration of priority-aware sensing and robust optimization remains largely unexplored. This motivates the incorporation of a priority map that assigns different importance levels to different regions, enabling sensor deployment strategies to preferentially safeguard critical infrastructures while maintaining robust coverage under uncertainty.

In this article, we propose a robust priority-aware coverage optimization framework for aerial sensor networks operating under sensor location uncertainty. Unlike conventional coverage optimization approaches that implicitly assume uniform importance across the monitored region, the proposed framework incorporates a priority map that assigns different importance levels to different regions, thereby enabling sensing resources to be preferentially allocated to critical areas. A mathematical optimization model is developed by integrating the weighted coverage objective with surveillance constraints and an RRF-based robustness formulation under sensor location uncertainty. To improve robustness against deployment inaccuracies, the radius of robust feasibility (RRF) \cite{rrf_1} is incorporated into the orientation optimization process. A priority-aware robust orientation optimization (PAROO) algorithm is then proposed to determine sensor orientations that maximize the weighted coverage objective. Experimental results on an airport-inspired surveillance scenario demonstrate that the proposed framework effectively directs sensing resources toward high-priority regions and achieves higher weighted coverage than representative baseline approaches, highlighting its potential for priority-aware aerial surveillance applications.

The remainder of this paper is organized as follows. Section~2 introduces the mathematical preliminaries, including the RRF, Minkowski operations, and the ASN model. The nominal priority-aware coverage optimization problem are presented in section~3. Section~4 develops the uncertainty model and derives the corresponding robust optimization framework, followed by the proposed priority-aware solution methodology. The experimental setup, simulation results, and comprehensive performance analysis of the proposed approach are reported in section~5. Finally, Section~6 concludes the paper and outlines potential directions for future research.

\section{Preliminaries}
Throughout this paper, $\mathbb{R}^n$ denotes the $n$-dimensional Euclidean space, while $\|\cdot\|$ represents the Euclidean norm and $\mathbb{B}_n$ denotes the open unit ball in $\mathbb{R}^n$. For any vector $a\in\mathbb{R}^n$, its transpose is written as $a^T$. Moreover, the Euclidean distance between two points $x,y\in\mathbb{R}^n$ is defined by
\[
d(x,y)=\|x-y\|.
\]
\subsection{Radius of robust feasibility}
Consider the following uncertain parametric linear system, represented by $\sigma^{\alpha}$:
\begin{equation}
    \sigma ^{\alpha } := \left\{ a^T_ix \leq b_i \right\}; \;(a_i,b_i) \in \mathcal{U} _i ^{\alpha},
\end{equation}
where $(a_i,b_i)$, $i=1,\cdots , m$, denote uncertain constraint parameters, and each uncertainty set $\mathcal{U}_i^{\alpha}$ is assumed to be a compact and convex subset of $\mathbb{R}^{n+1}$.

The robust counterpart (RC) is constructed by requiring every constraint to be satisfied for all admissible realizations of the uncertain parameters within the prescribed uncertainty sets. Consequently, when $\mathcal{U} _i^{\alpha}=(\bar{a}_i,\bar{b}_i)+\alpha \mathbb{B}_{n+1}$ for $i=1,\cdots , m$, the deterministic RC of $\sigma^{\alpha}$, together with its corresponding robust feasible region $F_R^{\alpha}$, is given by
\begin{equation}
    \sigma _ R ^ {\alpha} := \left\{ a_i^Tx \leq b_i ,\;  \forall (a_i,b_i) \in \mathcal{U} _i^{\alpha},\; i=1:m \right\}.
\end{equation}

Unlike deterministic optimization, which relies solely on nominal data, robust optimization explicitly accounts for uncertainty during the formulation stage. As a result, the computed solutions remain feasible for every realization contained within the specified uncertainty sets, thereby improving the reliability of the optimization model.

The radius of robust feasibility (RRF) characterizes the maximum allowable magnitude of uncertainty that preserves the feasibility of a given solution. It is formally defined as follows:
\begin{definition}(RRF)
Let $\sigma^{\alpha}$ be an uncertain parametric linear system and let $F_R^{\alpha}$ denote the feasible set corresponding to its RC. RRF is defined as the supremum of all uncertainty levels that preserve the feasibility of the robust system and is expressed as follows:
\begin{equation}
     \rho = \sup \{ \alpha \in \mathbb{R}_+ : (F_R ^ \alpha) \text{ is nonempty} \}.
\end{equation}
\end{definition}
The Minkowski function provides a convenient measure for scaling convex uncertainty sets and is defined as follows:

\begin{definition}(Minkowski function)
Let $\omega\subset\mathbb{R}^n$ be a convex set containing the origin in its interior. The associated Minkowski (gauge) function $\phi_{\omega}:\mathbb{R}^n\rightarrow\mathbb{R}_+$, where $\mathbb{R}_+=[0,\infty)$, is defined by:
    \begin{equation}
        \phi_{\omega} (x) := \inf{ \left\{ t>0 : x \in t \omega \right\} },\; x \in \mathbb{R} ^n.
    \end{equation}
\end{definition}
The RRF measures the largest allowable perturbation that preserves feasibility. Its exact expression for an arbitrary set $S\subset\mathbb{R}^{n+1}$ is obtained through the hypographical set $H(\bar{a},\bar{b})$ of the nominal system $\sigma$:
\begin{equation}
    H(\bar a, \bar b) := conv \left\{ (- \bar{a_i} , -\bar{b_i}): i=1:m \right\} + \mathbb{R} _+ (0_n, -1),
\end{equation}
where $\bar a=(\bar a_1,\bar a_2,\ldots,\bar a_m)\in(\mathbb{R}^n)^m$ and $\bar b=(\bar b_1,\bar b_2,\ldots,\bar b_m)\in\mathbb{R}^m$ denote the nominal constraint parameters. Provided that the nominal system is feasible, the corresponding RRF of $\sigma^{\alpha}$ is expressed as
\begin{equation}
    \rho = \inf_{(a,b) \in H(\bar{a},\bar{b})} \phi_Z (a,b),
\end{equation}
where $Z\subset\mathbb{R}^{n}$ is a compact convex set containing $0_n$ in its interior \cite{rrf_1}.

\subsection{Sensing model}
\begin{figure}
    \centering
    \includegraphics[width=0.75\textwidth]{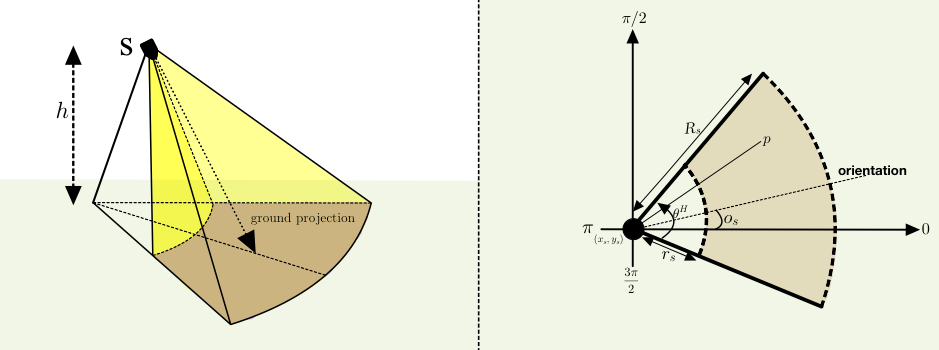}
    \caption{Illustration of the aerial directional sensing model and its projected sensing region on the ground plane.}
    \label{Fig_1}
\end{figure}
Let $\mathcal{X}\subset\mathbb{R}^{2}$ denote the monitoring region and $\mathcal{S}$ the set of deployed directional sensors. Each ground projection of sensor $s=(x_s,y_s)\in\mathbb{R}^{2}$ is characterized by its minimum and maximum sensing radii, $r_s$ and $R_s$, horizontal field of view $\theta^H$, and orientation angle $o_s$ measured from the positive $x$-axis (see Figure~\ref{Fig_1}). Accordingly, the sensing region is modeled as an annular sector, and a point $p=(x,y)\in\mathbb{R}^{2}$ is considered covered if it satisfies the corresponding radial and angular constraints. Thus, the effective sensing region $C_s(o_s)$ is defined as:
{\small
\begin{equation}
\label{sensing_area}
\resizebox{\linewidth}{!}{$
C_s(o_s)=
\left\{
p\in\mathbb{R}^{2}\,\;\middle|\,\;
r_s\leq ||s-p||\leq R_s,\,\,
\operatorname{atan2}(y-y_s,x-x_s)
\in
\left[
o_s-\tfrac{\theta^H}{2},\,
o_s+\tfrac{\theta^H}{2}
\right]
\right\}.
$}
\end{equation}
}
Accordingly, the sensing region of a directional sensor is completely specified by the parameters $(r_s,R_s,\theta^H,o_s)$.

\section{Problem formulation}
\begin{figure}
    \centering
    \includegraphics[width=\textwidth]{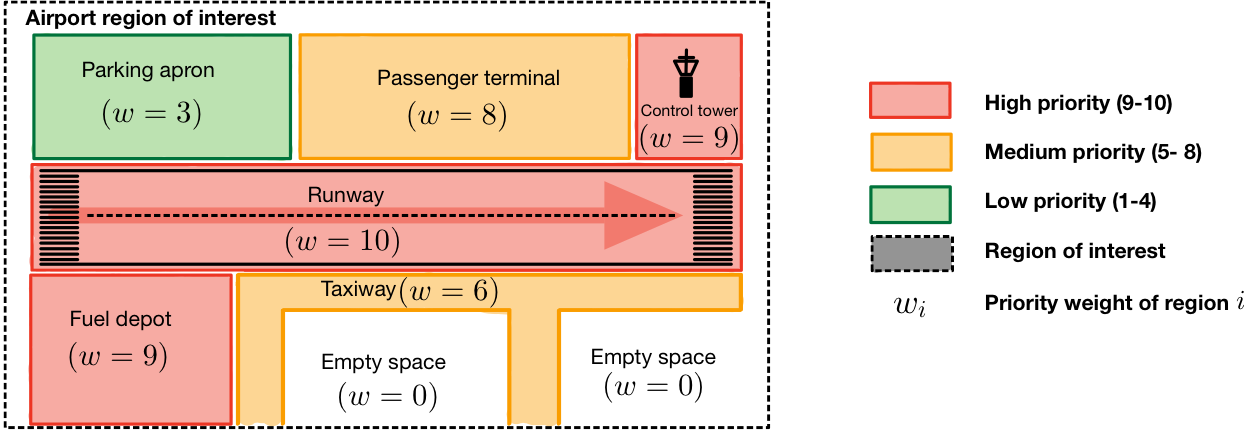}
    \caption{Airport surveillance scenario with priority-weighted monitoring regions used for coverage optimization.}
    \label{Fig_2}
\end{figure}
In practical surveillance applications, allocating sensing resources uniformly across the entire region of interest is often neither necessary nor efficient, since different locations possess different operational significance. Figure~\ref{Fig_2} illustrates a representative airport scenario, where facilities such as runways, fuel depots, and air traffic control towers require substantially higher monitoring reliability than terminals, taxiways, or parking areas. Therefore, each region is assigned a priority weight reflecting its relative importance, and the objective is to determine sensor orientations that maximize the weighted surveillance performance under sensor location uncertainty.

Let $\Omega \subset \mathbb{R}^{2}$ denote the region of interest (ROI) and let
$$\Omega=\bigcup_{k=1}^{K}\Omega_k,$$
where $\Omega_k$ represents the $k^{\mathrm{th}}$ priority region satisfying $\Omega_i\cap\Omega_j=\emptyset,\; i\neq j.$
Each region is assigned a positive priority weight according to its operational significance through the priority function
$\mathcal{W}:\Omega\rightarrow\mathbb{R}_{>0},$
defined as:
$$\mathcal{W} (p)=w_k,\qquad \forall p\in\Omega_k,$$
where $w_k>0$ denotes the priority weight associated with region $\Omega_k$. Larger values of $w_k$ correspond to regions requiring higher surveillance reliability, whereas smaller values indicate comparatively less critical areas.

Assume that $m$ aerial directional sensors are deployed over the ROI with their nominal location set as $\mathcal{S} = \{ s_1^o, s_2^o, \cdots, s_m^o\}$. The ground projected nominal position of the $i^{\mathrm{th}}$ sensor is denoted by
$s_i^o = (x_i^o,y_i^o)$, while its sensing orientation is represented by the corresponding pan angle $o_i$. Since the nominal sensor locations are assumed to be fixed during optimization, the orientation variables constitute the decision variables of the proposed optimization model. The objective is therefore to determine sensor orientations that maximize the weighted surveillance performance over the ROI.

Let $\mathbf{o}=(o_1,o_2,\ldots,o_m)$ denote the vector of sensor orientations, where $o_i$ is the orientation of sensor $s_i$. Moreover, let $c(t_j;\mathbf{o})\in\{0,1\}$ indicate whether point $t_j\in\Omega$ is covered under the orientation vector $\mathbf{o}$.
To ensure reliable sensing, the deployment must satisfy several operational constraints. First, a target (or grid point) is considered observable only if it lies within the effective sensing range of a sensor. Accordingly, let $\tau_{ij}$ denote the range feasibility indicator, defined by
\[
\tau_{ij}=
\begin{cases}
1,& r_i\le d(s_i^o,t_j)\le R_i,\\
0,& \text{otherwise},
\end{cases}
\]
where $r_i$ and $R_i$ denote the minimum and maximum sensing ranges of sensor $s_i^o$, respectively.

Furthermore, the target must lie inside the horizontal and vertical field of view of the sensor. Let $\theta ^H, \theta ^V$ denote the horizontal and vertical angles of view, respectively, and $\delta^H$ and $\delta^V$ denote the corresponding angular deviations between the sensor orientation and the target direction. The directional sensing requirements are given by
\begin{equation} \label{AOV}
    \begin{split}
        o_i - \frac{\theta ^H _i}{2} \leq \delta ^H _{ij} \leq o_i + \frac{\theta ^H _i}{2},\\
        o_i - \frac{\theta ^V _i}{2} \leq \delta ^V _{ij}  \leq o_i + \frac{\theta ^V _i}{2},
    \end{split}
\end{equation}
while the corresponding feasibility indicator is represented by
\[
\alpha_{ij}=
\begin{cases}
1,& \text{if Eq.~(\ref{AOV}) is satisfied},\\
0,& \text{otherwise}.
\end{cases}
\]

To guarantee effective surveillance quality, the inclination angle between the $i^{th}$ sensor and the $j^{th}$ target, say $\beta_{ij}$, should not exceed the prescribed threshold $\beta^{req}$. Accordingly,
\[
\beta_{ij}=
\begin{cases}
1,& \theta_{ij}\le\beta^{req},\\
0,& \text{otherwise}.
\end{cases}
\]

A target is regarded as successfully monitored only when all the above requirements are simultaneously satisfied. The nominal weighted coverage objective is then given by:
\begin{equation} \label{nominal}
    \text{(Nominal):}\;\max_{\mathbf{o}}\sum_{t_j\in\Omega}\mathcal{W}(t_j)\;c(t_j;\mathbf{o}).
\end{equation}
subject to the following surveillance constraint for all fixed $s_i^o$:
$$\gamma_{ij}=\tau_{ij} \times \alpha_{ij} \times \beta_{ij},
\qquad
\gamma_{ij}\in\{0,1\}.$$
Since the proposed formulation involves binary coverage variables and orientation optimization, it is NP-hard \cite{np_hard}. The resulting solution determines the sensor orientations that maximize target coverage under nominal conditions. However, the formulation assumes exact sensor locations. In practice, localization errors and environmental disturbances may shift the sensing footprints, degrading coverage performance. To ensure reliable deployment under such perturbations, sensor location uncertainty is incorporated into the robust optimization framework presented in the next section.

\section{Proposed methodology}

The nominal optimization model presented in the previous section is extended to account for sensor location uncertainty. Let the nominal position of sensor $s_i$ be denoted by $s_i^o=(x_i^o,y_i^o)$, while the actual deployed position is expressed as
\[
s_i=s_i^o+\Delta s_i,
\]
where $\Delta s_i$ is an unknown perturbation satisfying
$\|\Delta s_i\|\le\alpha$. The corresponding uncertainty set is defined by
\begin{align}
\label{uncertainity}
\mathcal U_i^\alpha=s_i^o+\alpha\mathbb B_2,
\end{align}
where $\alpha>0$ specifies the maximum admissible deviation from the nominal sensor location.

Accordingly, the nominal optimization model~\eqref{nominal} is transformed into the following robust counterpart:
\begin{equation}
\label{robust}
\text{(Robust counterpart):}\qquad
\max_{\mathbf o}
\sum_{t_j\in\Omega}
\mathcal W(t_j)c(t_j;\mathbf o),
\end{equation}
subject to
\begin{align}
\gamma_{ij}
=
\tau_{ij} \times \alpha_{ij} \times \beta_{ij},
\qquad
\gamma_{ij}\in\{0,1\},
\end{align}
for every admissible realization
$s_i\in\mathcal U_i^\alpha$.

The RRF introduced in Section~2 is then employed to determine the largest admissible perturbation for which all sensing constraints remain satisfied. Since evaluating infinitely many perturbations within $\mathcal U_i^\alpha$ is computationally impractical, the uncertainty region is approximated using a finite collection of boundary realizations generated according to the corresponding RRF.

Based on the computed RRF, the representative perturbed location of the $i^{\mathrm{th}}$ sensor along the sensing direction $\vec u$ is given by
\begin{equation}
s_i^\rho
=
s_i^o
+
\rho_i
\frac{\vec u}{\|\vec u\|},
\end{equation}
which provides a tractable approximation for evaluating sensing performance under admissible perturbations.

Let
\[
\mathcal B_i
=
\{b_1,b_2,\ldots,b_{N_i}\}
\]
denote the set of sampled boundary realizations associated with sensor $s_i$. For each candidate orientation, the weighted coverage objective is evaluated at every boundary sample, and the average value over all samples is used to assess its robustness. Consequently, the proposed optimization model is formulated as
\begin{equation}
\label{RC}
\text{(Proposed model):}\qquad
\max_{\mathbf o}
\frac1{|\mathcal B_i|}
\sum_{b_k\in\mathcal B_i}
\sum_{t_j\in\Omega}
\mathcal W(t_j)
c(t_j;b_k,\mathbf o),
\end{equation}
subject to
\begin{align}
\gamma_{ij}
=
\tau_{ij} \times \alpha_{ij} \times \beta_{ij},
\qquad
\gamma_{ij}\in\{0,1\},
\end{align}
for every
$b_k\in\mathcal B_i$.

To efficiently solve the proposed optimization problem, an iterative orientation optimization strategy is adopted. The monitoring region is discretized into uniform grid points, and the weighted coverage objective defined in Section~3 is evaluated over this discretization. Beginning from an initial deployment, the orientation of each sensor is updated sequentially while keeping the remaining sensors fixed. For every candidate orientation, the average weighted coverage over all sampled boundary realizations is computed, and the orientation yielding the largest objective value is retained. The procedure is repeated until no further improvement in the objective value is observed or a prescribed maximum number of iterations is reached.

The resulting optimization problem is solved using the proposed priority-aware robust orientation optimization (PAROO) algorithm, which iteratively determines robust sensor orientations by maximizing the average weighted coverage over the RRF boundary samples. The complete computational procedure is summarized in Algorithm~\ref{alg:PAROO}.

\begin{algorithm}
\caption{Priority-aware robust orientation optimization (PAROO)}
\label{alg:PAROO}
\begin{algorithmic}[1]

\Require Sensor set $\mathcal{S}$, monitoring grid $\Omega$, priority function $\mathcal{W}(\cdot)$, RRF values $\{\rho_i\}$, orientation step $\Delta o$, boundary sampling step $\Delta s$
\Ensure Optimized orientation vector $\mathbf{o}^{\ast}$

\State Initialize the orientation vector $\mathbf{o}$.

\For{each sensor $s_i\in\mathcal{S}$}

    \State $o_i^{\ast}\gets o_i,\qquad
    F_i^{\ast}\gets -\infty$

    \For{$o=0,\Delta o,\ldots,360-\Delta o$}

        \State $\overline{F}\gets0$
        \State Generate boundary samples $\mathcal{B}_i$ using $\rho_i$ and $\Delta s$.

        \For{each $b\in\mathcal{B}_i$}

            \State Evaluate $c(t_j;b,o)$ satisfying the sensing constraints.
            \State Compute
            \[
            F=\sum_{t_j\in\Omega}\mathcal W(t_j)c(t_j;b,o).
            \]
            \State $\overline{F}\gets\overline{F}+F$

        \EndFor

        \State $\overline{F}\gets\overline{F}/|\mathcal B_i|$

        \If{$\overline{F}>F_i^{\ast}$}

            \State $F_i^{\ast}\gets\overline{F}$
            \State $o_i^{\ast}\gets o$

        \EndIf

    \EndFor

    \State $o_i\gets o_i^{\ast}$

\EndFor

\State \Return $\mathbf{o}^{\ast}$

\end{algorithmic}
\end{algorithm}

\section{Experimental results and analysis}
The proposed framework is implemented in Python~3.12 using the Spyder IDE. Simulations are conducted over a two-dimensional monitoring region with randomly deployed aerial directional sensors. An airport-inspired priority map is adopted, where the runway, terminal, taxiway, and parking area are assigned different surveillance priorities according to their operational importance. The spatial distribution of these priorities is generated based on the airport layout shown in Figure~\ref{Fig_2}. Sensor location uncertainty is incorporated through the proposed RRF-based uncertainty model, and the PAROO algorithm is employed to maximize the weighted coverage objective.

\begin{figure}
\centering

\begin{subfigure}{0.48\textwidth}
    \centering
    \includegraphics[width=\linewidth]{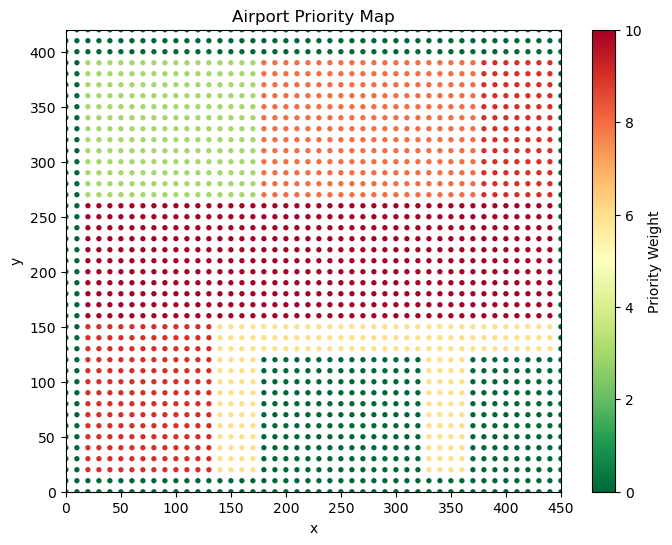}
    \caption{Airport priority map.}
\end{subfigure}

\begin{subfigure}{0.48\textwidth}
    \centering
    \includegraphics[width=\linewidth]{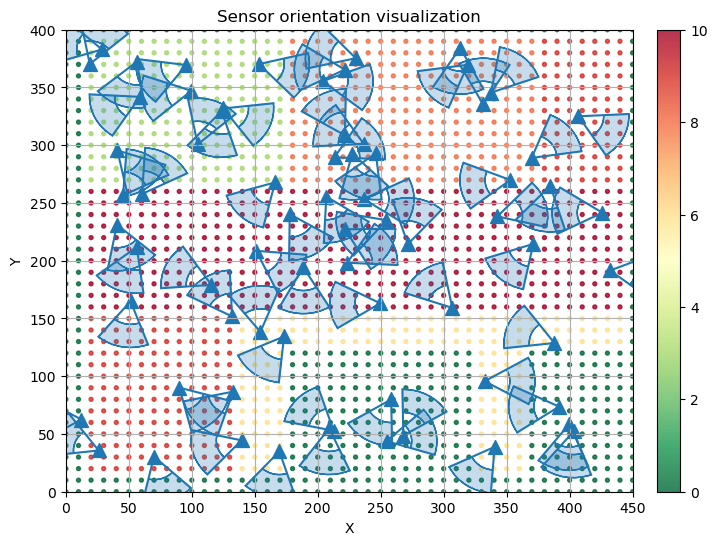}
    \caption{Initial deployment.}
\end{subfigure}
\hfill
\begin{subfigure}{0.48\textwidth}
    \centering
    \includegraphics[width=\linewidth]{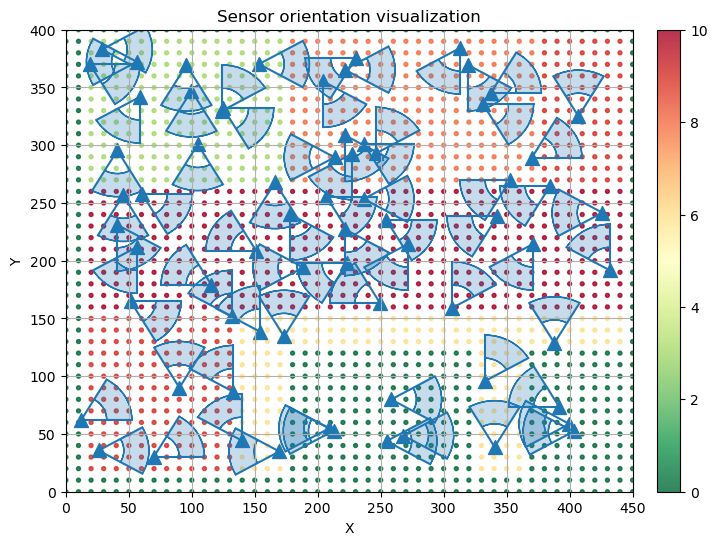}
    \caption{Deployment after PAROO.}
\end{subfigure}

\caption{Visualization of the proposed surveillance framework. (a) Airport-inspired priority map. (b) Initial sensor deployment with randomly oriented sensing sectors. (c) Optimized deployment obtained using the proposed PAROO algorithm, where sensing sectors are redirected toward high-priority airport facilities.}
\label{Fig_deployment}

\end{figure}
Figure~\ref{Fig_deployment} illustrates the generated priority map together with the sensor deployments before and after applying PAROO. It can be observed that the optimized sensing sectors are naturally redirected toward regions of higher surveillance priority while maintaining effective coverage under location uncertainty.

\begin{figure}
    \centering
    \includegraphics[width=0.75\textwidth]{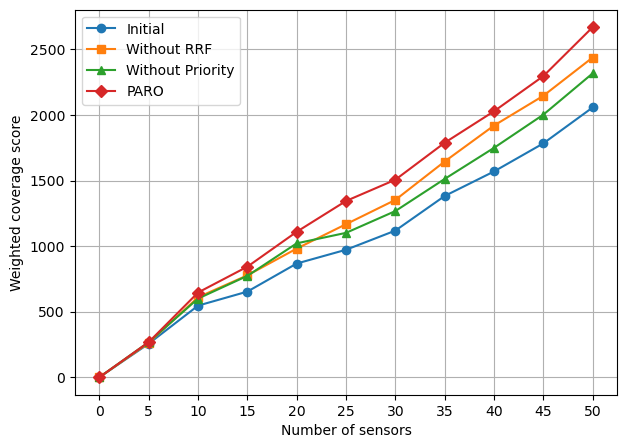}
    \caption{Comparison of the weighted coverage achieved by the initial deployment, the orientation optimization without RRF, the RO without priority-aware weighting, and the proposed PAROO algorithm.}
    \label{Fig_comparison}
\end{figure}
To assess the effectiveness of the proposed framework, PAROO is compared with three representative configurations: (i) the initial deployment, (ii) orientation optimization without the RRF-based robustness mechanism, and (iii) robust orientation optimization without priority-aware weighting. The corresponding weighted coverage obtained with increasing numbers of deployed sensors is illustrated in Figure~\ref{Fig_comparison}. The overall percentage improvements achieved by the different approaches are summarized in Remark~\ref{airport_remark}.

\begin{remark}{\rm
PAROO improves the weighted coverage by approximately $35.63\%$ over the initial deployment, compared with improvements of about $29.08\%$ and $17.75\%$ obtained without RRF and without priority-aware weighting, respectively. These results confirm the effectiveness of integrating RRF and priority-aware optimization.}\label{airport_remark}
\end{remark}

\section{Conclusion and future directions}
This paper proposed a priority-aware robust coverage optimization framework for ASNs under sensor location uncertainty. By integrating a weighted priority map into the coverage objective, the proposed framework allocates sensing resources according to the relative importance of different monitoring regions. To improve robustness against deployment inaccuracies, the RRF was incorporated into the orientation optimization process, and a priority-aware robust orientation optimization algorithm was developed to determine sensor orientations that maximize the average weighted coverage over admissible perturbations. Experimental results on an airport-inspired surveillance scenario demonstrated that the proposed framework effectively directs sensing sectors toward critical facilities and achieves higher weighted coverage than the considered baseline approaches.

Future work will focus on extending the proposed framework to dynamic priority maps, heterogeneous aerial sensor networks, and three-dimensional surveillance environments. Another promising direction is the development of distributed optimization algorithms for large-scale deployments with time-varying uncertainties.

\end{document}